\documentclass[a4paper,12pt]{article}
\usepackage{amsmath,amsfonts,amssymb,amsthm}

\begin{document}

\newcommand{\R}{\ensuremath{\mathbb{R}}}
\newcommand{\N}{\ensuremath{\mathbb{N}}}
\newcommand{\Z}{\ensuremath{\mathbb{Z}}}

\newcommand{\ep}{\varepsilon}
\newcommand{\lam}{\lambda}
\newcommand{\de}{\delta}
\newcommand{\De}{\Delta}

\newcommand{\LL}{{\cal L}}
\newcommand{\ST}{{\cal S}}
\newcommand{\UU}{{\cal U}}
\newcommand{\NN}{{\cal N}}
\newcommand{\GG}{{\cal G}}
\newcommand{\WW}{{\cal W}(\delta_1,\delta_2,p,p')}

\newcommand{\ds}{\mbox{dist}}

\newcommand{\pp}{P(\delta_1,\delta_2,p)}
\newcommand{\qp}{Q(\delta_1,\delta_2,p)}
\newcommand{\tp}{T(\delta_1,\delta_2,p)}
\newcommand{\ttp}{T(\delta_1,\Delta,p)}
\newcommand{\ppd}{P(\delta_1,\Delta,p)}
\newcommand{\spp}{S(\delta_1,\Delta,p)}
\newcommand{\ip}{\mbox{Int}^0\;}
\newcommand{\bnp}{\partial^0 P}
\newcommand{\bnq}{\partial^0 Q}

\newcommand{\ppf}{P(\delta_1,\delta_2,f(p))}
\newcommand{\qpf}{Q(\delta_1,\delta_2,f(p))}

\newcommand{\ppp}{P(\delta_1,\delta_2,p')}
\newcommand{\qpp}{Q(\delta_1,\delta_2,p')}

\renewcommand{\thefootnote}{\fnsymbol{footnote}}

\title{Lyapunov functions, shadowing, and topological stability}

\author{Alexey A.\ Petrov\footnotemark[1] \hspace{0.15cm}\footnotemark[2] 
\hspace{0.15cm}, Sergei Yu.\ Pilyugin\footnotemark[1]
\hspace{0.15cm}\footnotemark[2]}
\date{\today}

\footnotetext[1] {Faculty of Mathematics and Mechanics, St.\ Petersburg State 
University, University av.,\ 28, 198504, St.\ Petersburg, Russia,
e-mail: sp@sp1196.spb.edu, al.petrov239@gmail.com.}

\footnotetext[2] {Supported by the Russian Foundation for Basic Research
(project 12-01-00275); the first author is also supported by the Chebyshev
Laboratory, Faculty of Mathematics and Mechanics, St.\ Petersburg State 
University, under grant 11.G34.31.0026 of the Government of Russian Federation.}
\maketitle

\renewcommand{\thefootnote}{\arabic{footnote}}

\begin{abstract}
\center{We use Lyapunov type functions to give new conditions
under which a homeomorphism of a compact metric space has the
shadowing property. These conditions are applied to establish
the topological stability of some homeomorphisms with nonhyperbolic
behavior. Bibliography: 10 titles.}
\end{abstract}


{\bf AMS(MOS) subject classifications.}  54C60, 37C75\\

\section{Introduction}

The shadowing property of dynamical systems (diffeomorphisms or flows)
is now well-studied (see, for example, the monographs [4, 6] and the
recent survey [7]).
This property means that, near approximate trajectories
(so-called pseudotrajectories), there exist exact trajectories of the
system.

Mostly, standard methods allow one to show that the shadowing property
follows from hyperbolic behavior of trajectories of the system.
It is well known that a structurally stable system has the
shadowing property (and this property is Lipschitz), see [6].

One can mention several papers which contain methods of proving
the shadowing property for systems with nonhyperbolic behavior
(see, for example, [1]).

In [3], Lewowicz used Lyapunov type functions to study topological
stability of dynamical systems (see also [8]). This property is stronger
than the shadowing property (and they are equivalent for expansive
systems on smooth closed manifolds, see Sec.~3).

In this paper, we give sufficient conditions of shadowing for a
dynamical system generated by a homeomorphism of a compact metric space in terms
of existence of a pair of Lyapunov type functions.
These conditions are formulated in terms of some sets related to the
considered pair of functions.

In fact, our conditions have much in common with the topological
conditions used in the classical Wazewski principle in the theory of
differential equations [10]. In a sense, a close reasoning has been
used by the second author in his joint paper [5] with Plamenevskaya
devoted to the $C^0$-genericity of shadowing.

The structure of the paper is as follows. In Sec.~2, we formulate and
prove our basic shadowing result. Section~3 is devoted to topological
stability. We consider an example studied by Lewowicz in [3] and
give a comment on the method of [3]. One more example shows that
our methods are applicable to homeomorphisms.

\section{Lyapunov functions and shadowing}

Let $f$ be a homeomorphism of a metric space $(X,\ds)$. 

As usual, we say that a sequence $\{p_k\in X:\;k\in\Z\}$ is a $d$-pseudotrajectory
of $f$ if
\begin{equation}
\label{s00}
\ds(p_{k+1},f(p_k))<d,\quad k\in\Z.
\end{equation}
We say that a pseudotrajectory $\{p_k:\;k\in\Z\}$ is $\ep$-shadowed by
a point $r$ if
$$
\ds(f^k(r),p_k)<\ep, \quad k\in\Z.
$$
We say that $f$ has the standard shadowing property if for any $\ep>0$
we can find a $d>0$ such that any $d$-pseudotrajectory of $f$ is
$\ep$-shadowed by some point.

It is well known (see [6]) that to establish the standard shadowing property
on a compact phase space
it is enough to show that $f$ has the so-called finite shadowing property:
For any $\ep>0$ we can find a $d>0$ (depending on $\ep$ only) such that
if $\{p_k:\;0\leq k\leq m\}$ is a finite $d$-pseudotrajectory, then there 
is a point $r$ such that
\begin{equation}
\label{s01}
\ds(f^k(r),p_k)<\ep, \quad 0\leq k\leq m.
\end{equation}

Our goal is to give sufficient conditions under which a homeomorphism
has the finite shadowing property on $X$.
In our conditions, we use analogs of Lyapunov functions.

Let us formulate our main assumptions.

We assume that the space $X$ is compact and there exist two continuous
nonnegative functions $V$ and $W$ defined in a closed neighborhood of the
diagonal of $X\times X$ such that $V(p,p)=W(p,p)=0$ for any $p\in X$
and the conditions (C1)-(C9) stated below are satisfied.
In what follows, arguments of the functions $V$ and $W$ are assumed
to be close enough, so that the functions are defined.

We formulate our conditions not directly in terms of the
functions $W$ and $V$ but in terms of some geometric objects
defined via these functions. Our main reasoning for the choice of this
form of conditions is as follows:

(1) Precisely these conditions are used in the proofs;

(2) it is easy to check conditions of that kind for particular
functions $W$ and $V$ (see the examples below).

Fix positive numbers $a,b>0$ and a point $p\in X$ and let
$$
P(a,b,p)=\{q\in X:\;V(q,p)\leq a,\;W(q,p)\leq b\},
$$
$$
Q(a,b,p)=\{q\in P(a,b,p):\;\;V(q,p)=a\},
$$
and
$$
T(a,b,p)=\{q\in P(a,b,p):\;V(q,p)=0\}.
$$

Denote by $B(\ep,p)$ the open $\ep$-ball centered at $p$.

Set
$$
\mbox{Int}^0\,P(a,b,p)=\{q\in P(a,b,p):\;V(q,p)<a,\;W(q,p)<b\},
$$
$$
\bnp(a,b,p)=Q(a,b,p)\cup \{q\in P(a,b,p):\;W(q,p)=b\},
$$
and
$$
\ip Q(a,b,p)=\{q\in P(a,b,p):\;V(q,p)=a,\;W(q,p)<b\}.
$$

Conditions (C1) -- (C4) contain our assupmtions on the
geometry of the sets introduced above.

(C1) For any $\ep>0$ there exists a $\Delta_0=\De_0(\ep)>0$ such that 
$P(\De_0,\De_0,p)\subset B(\ep,p)$ for $p\in X$.

There exists a $\De_1>0$ such that if $p\in X$,
$\de_1,\de_2,\De<\De_1$, and $\de_2<\De$,
then there exists a number $\alpha=\alpha(\de_1,\de_2,\De)>0$ such that

(C2) $\qp$ is not a retract of  $\pp$;

(C3) $\qp$ is a retract of  
$\pp\setminus\tp$;

(C4) there exists a retraction
$$
\sigma:\ppd\to \pp
$$
such that $V(\sigma(q),p)\geq \alpha V(q,p)$ for $q\in\ppd$.

In the next group of conditions, we state our assumptions on
the behavior of the introduced objects and their images
under the homeomorphism $f$.

We assume that for any $\De<\De_1$ there exists positive numbers 
$\de_1,\de_2<\De$ such the following relations hold for
any $p\in X$:

(C5) $f(\pp)\subset \ip P(\De,\De,f(p))$,

$f^{-1}(\ppf)\subset \ip P(\De,\De,p)$;

(C6) $f(\tp)\subset \ip\,P(\de_1,\de_2,f(p))$;

(C7) $f(\ttp)\cap\qpf=\emptyset$;

(C8) $f(\pp)\cap \bnp(\de_1,\de_2,f(p))\subset\ip Q(\de_1,\de_2,f(p))$;

(C9) $f(\spp)\cap\ppf=\emptyset$, where

$\spp=\{q\in P(\De,\De,p)\,:V(q,p)\geq \de_1$\}.
\medskip

Our main result is as follows.
\medskip

{\bf Theorem 1. }{\em Under conditions} (C1)-(C9), $f$ {\em has the
finite shadowing property on the space $X$}.
\medskip

In the proof of this statement, we apply the following two lemmas.

First let us formulate one more condition (the letter ${\cal W}$
in this condition indicates that, as was mentioned above, this condition
has much in common with the classical Wazewski principle in
the theory of differential equations).
\medskip

Let $p,p'\in X$ and $\de_1,\de_2>0$. We say that
condition $\WW$ holds if
\begin{equation}
\label{s4}
f(P)\cap\partial^0 P'\subset Q',
\end{equation}
\begin{equation}
\label{s5}
f(Q)\cap P'=\emptyset,
\end{equation}
and $Q$ is a retract of the set 
$$
H=H_1\cup f^{-1}(Q'),
$$ 
where $P=\pp$, $Q=\qp$, $P'=\ppp$, $Q'=\qpp$, and $H_1=P\setminus
f^{-1}(\ip P')$.
\medskip

{\bf Lemma~1. }{\em Let positive numbers $\de_1,\de_2<\De$ satisfy
conditions} (C4) -- (C9). {\em Let $\de=\min(\de_1,\de_2)$.
There exists
a positive $d=d(\delta)$ such that
if} $\ds(p',f(p))<d$, {\em then
condition $\WW$ holds}.
\medskip

{\em Proof. } In the proof, we several times select a small $d$ (depending
on $\de$) and then take as the required $d$ the minimum of the selected values of $d$.

Condition (C6), the compactness of the neighborhood of the diagonal of $X\times X$
in which the functions $V$ and $W$ are defined, and the continuity of $f$
imply that there exist positive numbers $c_1<\de_1$ 
and $c_2<\de_2$ such that
if $q\in f(\tp)$, then $V(q,f(p))\leq c_1$ and $W(q,f(p))\leq c_2$.
Hence, there exists a $d=d(\delta)$ such that if
\begin{equation}
\label{s8}
\ds(p',f(p))<d,
\end{equation}
then $V(q,p')<\de_1$ and $W(q,p')<\de_2$, which means that
\begin{equation}
\label{s9}
f(\tp)\subset \mbox{Int}^0\,P'.
\end{equation}

A similar reasoning based on condition (C5) shows that 
there exists a $d=d(\delta)$ such that if inequality (\ref{s8})
is satisfied, then
\begin{equation}
\label{s7}
f(P)\subset P(\De,\De,p')
\end{equation}
and
\begin{equation}
\label{s10}
f^{-1}(P')\subset P(\De,\De,p).
\end{equation}

In particular, inclusion (\ref{s10}) implies that
\begin{equation}
\label{s11}
f^{-1}(Q')\subset P(\De,\De,p).
\end{equation}
Let us show that 
there exists a $d=d(\delta)$ such that if inequality (\ref{s8})
is satisfied, then
\begin{equation}
\label{s12}
f^{-1}(Q')\subset P(\de_1,\De,p).
\end{equation}

Since the set $S:=\spp$ is compact, it follows from condition (C9)
that there exists a number $c_3>0$ such that
if $q\in S$, then 
$$
\max(V(f(q),f(p))-\de_1,W(f(q),f(p))-\de_2)\geq c_3.
$$
Hence, 
there exists a $d=d(\delta)$ such that if inequality (\ref{s8})
is satisfied, then 
$$
\max(V(f(q),p')-\de_1,W(f(q),p')-\de_2)>0
$$
for $q\in S$, which 
implies that condition (\ref{s5}) is satisfied and
$f(S)\cap Q'=\emptyset$. Now inclusion (\ref{s12})
follows from inclusion (\ref{s11}).

Clearly, condition (C8) (combined with inclusion (\ref{s7}))
implies that there 
exists a $d=d(\delta)$ such that if inequality (\ref{s8})
is satisfied, then inclusion (\ref{s4}) holds.

Similarly, it follows from condition (C7) that 
there exists a number $c_4>0$ such that
if $q\in f^{-1}(Q)$, then $V(q,p)\geq 2c_4$. Hence, 
there exists a $d=d(\delta)$ such that if inequality (\ref{s8})
is satisfied, then $V(q,p)\geq c_4$ for $q\in f^{-1}(Q')$.

Apply condition (C4) to find a retraction
$$
\sigma:\,\ppd\to \pp
$$
such that 
\begin{equation}
\label{s14}
V(\sigma(q),p)\geq \alpha c_3,\quad q\in f^{-1}(Q').
\end{equation}
The set
$$
U=\sigma(f^{-1}(Q'))
$$
is compact, and inequality (\ref{s14}) implies that
\begin{equation}
\label{s13}
U\cap \ttp=\emptyset.
\end{equation}

By condition (C3), there exists a retraction $\rho_0$
of $P\setminus \tp$ to $Q$. Relations (\ref{s9}) and (\ref{s13})
imply that $H_1\cup U\subset P\setminus \tp$.
Hence, the restriction of
$\rho=\rho_0\circ\sigma$ to $H$ is the required
retraction $H\to Q$. $\bullet$
\medskip

{\bf Lemma~2. } {\em Let $p_0,\dots,p_m$ be points in $X$ such that
condition ${\cal W}(\de_1,\de_2,p_k,p_{k+1})$ holds for $k=0,\dots,m-1$.
Then there exists a point $r\in P(\de_1,\de_2,p_0)$ such that
$f^k(r)\in P(\de_1,\de_2,p_k)$ for $k=1,\dots,m$.}
\medskip

{\em Proof. } Consider the sets
$$
A_k=P(\de_1,\de_2,p_k)\setminus\bigcap_{l=k+1}^m f^{-(l-k)}
(\mbox{Int}^0\,P(\de_1,\de_2,p_l)),\quad k=0,\dots, m-1.
$$
It follows from equality (\ref{s5}) that
$$
f(Q(\de_1,\de_2,p_k))\cap P(\de_1,\de_2,p_{k+1})=\emptyset.
$$
Hence, $Q(\de_1,\de_2,p_k)\subset A_k$.

We claim that there exist retractions
$$
\rho_k:\,A_k\to Q(\de_1,\de_2,p_k),\quad k=0,\dots,m-1.
$$
This is enough to prove our lemma since the existence of 
$\rho_0$ means that
$$
\bigcap_{l=0}^m f^{-l}(\mbox{Int}^0\,P(\de_1,\de_2,p_l))\neq\emptyset
$$
(otherwise there exists a retraction of $P(\de_1,\de_2,p_0)$ to
$Q(\de_1,\de_2,p_0)$, which is impossible by condition (C3)).

The existence of $\rho_{m-1}$ is obvious since condition
${\cal W}(\de_1,\de_2,p_{m-1},p_{m})$
implies the existence of a retraction 
$$
A_{m-1}\cup f^{-1}(Q(\de_1,\de_2,p_m))\to Q(\de_1,\de_2,p_{m-1}).
$$

Let us assume that the existence of retractions $\rho_{k+1},\dots,\rho_{m-1}$
has been proved. Let us prove the existence of $\rho_k$.

For brevity, we denote $P_k=P(\de_1,\de_2,p_k), Q_k=Q(\de_1,\de_2,p_k), 
P_{k+1}=P(\de_1,\de_2,p_{k+1})$,
and $Q_{k+1}=Q(\de_1,\de_2,p_{k+1})$.

Note that the definition of the sets $A_k$ implies that
\begin{equation}
\label{s15}
A_k\cap f^{-1}(P_{k+1})\subset f^{-1}(A_{k+1}).
\end{equation}

Define a mapping $\theta$ on $A_k$ by setting
$$
\theta(q)=f^{-1}\circ\rho_{k+1}\circ f(q),\quad
q\in A_k\cap f^{-1}(P_{k+1}),
$$
and
$$
\theta(q)=q, \quad q\in A_k\setminus f^{-1}(P_{k+1}).
$$
Inclusion (\ref{s15}) shows that the mapping $\theta$
is properly defined.

Let us show that this mapping is continuous. Clearly,
it is enough to show that $\rho_{k+1}(r)=r$ for
$r\in f(A_k\cap f^{-1}(\partial^0 P_{k+1}))$.
For this purpose, we note that
$$
f(A_k\cap f^{-1}(\partial^0 P_{k+1}))=f(A_k)\cap\partial^0 P_{k+1}
\subset f(P_k)\cap\partial^0 P_{k+1}\subset Q_{k+1}
$$
(we refer to inclusion (\ref{s4})) and $\rho_{k+1}(r)=r$
for $r\in Q_{k+1}$.

Clearly, $\theta$ maps $A_k$ to the set
\begin{equation}
\label{s16}
[P_k\setminus f^{-1}(P_{k+1})]\cup f^{-1}(Q_{k+1}).
\end{equation}

Condition ${\cal W}(\de_1,\de_2,p_k,p_{k+1})$ 
implies that there exists a retraction $\rho$ 
of (\ref{s16}) to $Q_k$. It remains to note that $\theta(q)=q$
for $q\in Q_k$ due to condition (\ref{s5}). Thus,
$$
\rho_k=\rho\theta:\,A_k\to Q_k
$$
is the required retraction. The lemma is proved. $\bullet$

To complete the proof of the main theorem, we take an arbitrary $\ep>0$,
apply condition (C1) to find a proper $\De_0$
and then find the corresponding numbers $\de_1,\de_2,\De$.
Lemma~1 implies that there exists a $d>0$ depending on $\de_1,\de_2,\De$
(i.e., on $\ep$) such that if $p_0,\dots,p_m$ is a finite
$d$-pseudotrajectory of $f$, then 
condition ${\cal W}(\de_1,\de_2,p_k,p_{k+1})$ holds for $k=0,\dots,m-1$.
Now it follows from Lemma~2 that $f$ has the finite shadowing
property on $X$. $\bullet$

\section{Topological stability}

In this section, we assume, for simplicity of presentation, that $X$ is a smooth
closed manifold. 

Let $H(X)$ be the space of homeomorphisms of $X$
endowed with the metric
$$
\rho(f,g)=\max_{p\in X}\max(\ds(f(p),g(p)),\ds(f^{-1}(p),g^{-1}(p))).
$$
It is well known that $H(X)$ is a complete metric space.

A homeomorphism $f$ is called topologically stable if for any
$\ep>0$ there exists a neighborhood $Y$ of $f$ in $X$ such that
if $g\in Y$, then there exists a continuous map $h:X\to X$
such that $f\circ h=h\circ g$ and
$$
\ds(h(p),p)<\ep,\quad p\in X.
$$

It is not difficult to show that if a homeomorphism $f$ is topologically
stable, then $f$ has the shadowing property (see [6]).

To formulate general sufficient conditions of topological stability,
let us recall one more standard definition.

A homeomorphism $f$ is called expansive if there exists a positive
number $a$ such that if
$$
\ds(f^k(p),f^k(q))\leq a,\quad k\in\Z,
$$
then $p=q$.

Walters proved in [9] the following theorem (see also [6]).
\medskip

{\bf Theorem 2. }{\em If a homeomorphism $f$ is expansive 
and has the shadowing property, then $f$ is topologically 
stable.}
\medskip

{\bf Remark.} Usually, an expansive homeomorphism having
the shadowing property is called topologically Anosov.
\medskip

It is known that structurally stable diffeomorphisms are
topologically stable [2].

We want to show that our shadowing result based on analogs
of Lyapunov functions can be applied in the proof of
topological stability of dynamical systems that are ``far"
from the set of structurally stable diffeomorphisms
(in particular, such systems may have nonhyperbolic fixed
points).
\medskip

{\bf Example 1. } Consider an example of a diffeomorphism
of the 2-torus $T^2$ studied by Lewowicz in [3].
This diffeomorphism is a perturbation of a hyperbolic automorphism of $T^2$.

Consider numbers $0<\alpha<1<\beta$ and a small $r>0$ and define
a map $F:\R^2\to\R^2$ by
$$
F(x,y)=(\alpha x+\lambda(x)\mu(y),\beta y),
$$
where
$$
\lambda(x)=\int_0^x((1-\alpha)-h(s))ds,
$$
$h:\R\to \R$ is a $C^\infty$ function such that
$h(0)=0$, $0\leq h(x)<1$, and 
$\lambda(x)=0$ for $|x|\geq r$;
$\mu:\R\to \R$ is a $C^\infty$ function such that
$\mu(0)=1$, $\mu(y)=\mu(-y)$, $\mu$ is not increasing 
for $y\geq 0$, and $\mu(y)=0$ for $|y|\geq r$.

Let $A$ be an integer hyperbolic $2\times 2$ matrix with $\mbox{det}A=1$.
If $0<\alpha<1<\beta$ are the eigenvalues of $A$
and $u_1$ and $u_2$ are the corresponding eigenvectors,
then
$$
A(x,y)=(\alpha x,\beta y)
$$
in cooordinates whose axes are parallel to $u_1$ and $u_2$.

The lattice $\Xi$ with vertices
$$
\{(n+1/2)u_1,(m+1/2)u_2\,:n,m\in\Z\}
$$
is invariant with respect to the action of the map $v\mapsto Av$.
Let $\pi:\R^2\to\R^2/\Xi$ be the corresponding projection of the plane to
the 2-torus.

Define
$f:T^2\to T^2$ by $f(\pi(\xi,\eta))=\pi\circ F(x,y)$
(of course, we extend $F$ periodically with respect to the
above-mentioned lattice).

It is shown in [3] that if $r$ is small enough, then $f$ is an
expansive diffeomorphism of the torus. At the same time, $f$ is
not Anosov (and is not structurally stable) since the eigenvalues
of $Df$ at the zero fixed point are $1$ and $\beta$.

Consider the functions $V$ and $W$ defined as follows.
If $p=(p_x,p_y)$ and $q=(q_x,q_y)$, we set $V(p,q)=|p_y-q_y|$
and $W(p,q)=|p_x-q_x|$ (such functions are properly defined
if $p$ and $q$ are close enough).
\medskip

It is obvious that conditions (C1) -- (C5) are satisfied.
Let us check condition (C6). We fix $0<\De<\De_0$ and a point 
$p=(p_x,p_y)$. Let $f=(f_x,f_y)$. Take a small $\de_2>0$
and consider points $p'=(p_x+\nu,p_y)$, where $|\nu|\leq\de_2$.
If $\nu=0$, then
$$
|f_x(p)-f_x(p')|=0.
$$
If $\nu\neq0$, then
$$
|f_x(p)-f_x(p')|=
|\nu(\alpha+\mu(p_y)(1-\alpha))-\mu(p_y)\int_{p_x}^{p_x+\nu}h(s)ds|<
$$
$$
<|\nu|\leq d_2
$$
(we take into account that $\mu(y)\leq 1$ and $h(x)>0$
for $x\neq 0$).

It follows that condition (C6) is satisfied for any 
$(\de_1,\de_2)$, where  $\de_2<\De$ is small enough.
In addition, if $\de_2$ is small enough, then the
``rectangle" $\pp$ is close to the ``segment" $\tp$,
which implies that condition (C8) is satisfied as well.

Since $f$ expands in the $y$ direction, conditions (C7)
and (C9) are satisfied automatically.

Thus, we can apply the main theorem to show that $f$
has the shadowing property, which implies the topological
stability of $f$.
\medskip

{\bf Remark. } Let us make a comment concerning the conditions
and proofs in the paper [3]. First, the proof in [3] refers to
the smoothness of the system considered (while our proof
works for homeomorphisms). Second, the proof in [3]
reduces the problem to study of suspension flows, which does 
not seem natural. Third, in our opinion, the proof in [3]
requires stronger assumptions on the regularity of the function $V$
than stated.
\medskip

As was mentioned, our methods are applicable to homeomorphisms.
\medskip

{\bf Example 2. } Consider a perturbation $f$ of the 
hyperbolic automorphism of $T^2$ corresponding to the map
$$
F(x,y)=(\mu_1(x),\mu_2(y)),
$$
where $\mu_1$ and $\mu_2$ are increasing continuous functions
for which there exist numbers $r,\lambda\in(0,1)$ such that

(1) 
$$
|\mu_1(x+\nu)-\mu_1(x)|\leq\lambda|\nu|
$$
and 
$$
\lambda^{-1}|\nu|\leq|\mu_2(y+\nu)-\mu_2(y)|
$$
for $|\nu|<r$;

(2) $\mu_1(x)=\alpha x,\quad |x|\geq r$;

(3) $\mu_2(y)=\beta y,\quad |y|\geq r$.

To prove that $f$ is topologically stable, one can apply the
same functions $V$ and $W$ and the same reasoning as in Example 1
(to show that $f$ is expansive, one can apply the same reasoning
as that applied in [3] to Example 1 considering the function
${\bf V}(p,q)=V(p,q)-W(p,q)$).

\section{References}

\noindent

1. S. M. Hammel, J. A. Yorke, and C. Grebogi,
{\em Numerical orbits of chaotic dynamical processes represent
true orbits}, Bull. Amer. Math. Soc. {\bf 19} (1988), 465-470.

2. M. Hurley, {\em Combined structural and topological stability are
equivalent to Axiom A and the strong transversality condition},
Ergod. Theory Dyn. Syst. {\bf 4} (1984), 81-88.

3. J. Lewowicz, {\em Lyapunov functions and topological stability},
J. Differential Equations {\bf 38} (1980), 192-209.

4. K. Palmer, {\em Shadowing in Dynamical Systems. Theory and
Applications}, Kluwer, Dordrecht, 2000.

5. S. Yu. Pilyugin and O. B. Plamenevskaya, {\em Shadowing is
generic}, Topology and Its Appl. {\bf 97} (1999), 253-266.

6. S. Yu. Pilyugin, {\em Shadowing in Dynamical Systems},
Lecture Notes in Math., vol. 1706, Springer, Berlin, 1999.

7. S. Yu. Pilyugin, {\em Theory of pseudo-orbit shadowing
in dynamical systems}, Differential Equations, {\bf 47}
(2011), 1929-1938.

8. J. Tolosa, {\em The method of Lyapunov functions in two
variables}, Contemp. Math. {\bf 440} (2007), 243-271.

9. P. Walters, {\em On the pseudo orbit tracing property and its
relationship to stability}, in: {\em The Structure of Attractors
in Dynamical Systems}, 
Lecture Notes in Math., vol. 668, Springer, Berlin, 1999, 231-244.

10. T. Wazewski, {\em Sur un principe topologique de l'examen
de l'allure asymptotique des int\'egrales des \'equations
differentielles ordinares}, Ann. Soc. Polon. Math., {\bf 20}
(1947), 279-313.

\end{document}